\documentclass[11pt,a4paper]{amsart}
\usepackage{amsmath}
\usepackage{amssymb}
\usepackage{amsfonts}
\usepackage{amsthm}
\usepackage{amsxtra}
\usepackage{amsopn}
\usepackage{array}
\usepackage{booktabs}
\usepackage{hyperref}
\hypersetup{
    colorlinks=true,
    linkcolor=blue,
    filecolor=magenta,
    urlcolor=blue
}

\newcommand{\fH}{\mathfrak{H}}
\newcommand{\cA}{\mathcal{A}}
\newcommand{\cH}{\mathcal{H}}
\newcommand{\C}{\mathbb{C}}
\newcommand{\R}{\mathbb{R}}

\DeclareMathOperator{\Log}{Log}
\DeclareMathOperator{\diag}{diag}

\theoremstyle{plain}
\newtheorem{thm}{Theorem}[section]
\newtheorem{lem}[thm]{Lemma}
\newtheorem{prop}[thm]{Proposition}
\newtheorem{cor}[thm]{Corollary}
\numberwithin{equation}{section}
\theoremstyle{definition}
\newtheorem{defn}[thm]{Definition}
\newtheorem{hyp}[thm]{Hypothesis}
\newtheorem{rem}[thm]{Remark}

\usepackage{fancyhdr}
\begin{document}
\title[Qc mappings between balls]{Quasiconformal mappings between Kor\'anyi balls, $CC$-balls and bubble sets in the first Heisenberg group}
\author{Ioannis D. Platis}

\address{Department of Mathematics,
University of Patras,
GR24508 Rion, Achaia, Greece.}
\email{idplatis@upatras.gr}


\begin{abstract}
We study quasiconformal mappings between fundamental model domains in the Heisenberg group $\fH$, including Kor\'anyi balls, Carnot-Carath\'eodory balls, and non-spherical bubble sets. Using a coordinate system based on horizontal lifts of simple plane curves, we foliate $\fH$ by surfaces of revolution. This geometric model allows us to study contact quasiconformal transformations known as spiral stretch-maps. We derive explicit formulas for their complex Beltrami coefficients $\mu_f$ and their maximal dilatations $K_f$, and compute their values along the equator of the profile curve. A universality emerges: for all pairwise combinations of the three domain families, the limiting Beltrami coefficient depends only on the stretching parameter $k$ and is given by $|k^2-1|/(k^2+1)$, with maximal dilatation $\max(k^2,1/k^2)$. Neither the domain type nor the target type affects the limit. In fact, the theorem holds for \emph{any} pair of profile curves satisfying a mild structural hypothesis stated in Section \ref{sec-profiles}. A brief numerical check for $k=2$ is included in Appendix \ref{app:num}.
\end{abstract}

\maketitle

\section{Introduction and Statement of Main Results}

\subsection{Sub-Riemannian geometry and the Heisenberg group}
Sub-Riemannian geometry studies smooth manifolds equipped with a bracket-generating distribution of tangent planes, together with a smoothly varying inner product on that distribution. The simplest and most fundamental example is the first Heisenberg group $\fH = \fH^1$, a step-two nilpotent Lie group diffeomorphic to $\R^3$ but equipped with a non-integrable horizontal distribution $\cH$ of rank two. The Heisenberg group serves as the model space for sub-Riemannian geometry in the same way that Euclidean space serves as the model for Riemannian geometry. Its geometry is genuinely different from Euclidean geometry: the Hausdorff dimension of $\fH$ with respect to the sub-Riemannian metric is four rather than three, the isoperimetric profile is non-Euclidean, and the group admits a one-parameter family of parabolic dilations $\delta_\lambda(z,t)=(\lambda z,\lambda^2 t)$ that reflect the distinct scaling of horizontal and vertical directions.

The horizontal distribution $\cH$ is the kernel of the contact form
$$
\omega = dt + 2(x\,dy - y\,dx) = dt + 2\Im(\overline z\,dz),
$$
and is spanned by the left-invariant vector fields
$$
X = \partial_x + 2y\,\partial_t, \qquad Y = \partial_y - 2x\,\partial_t.
$$
The bracket $[X,Y] = -4\partial_t$ generates the missing vertical direction, so that $\cH$ is bracket-generating. A curve $\gamma:[a,b]\to\fH$ is called \emph{horizontal} if $\dot\gamma(s)\in\cH_{\gamma(s)}$ for almost every $s$; equivalently, $\omega(\dot\gamma)=0$, which in coordinates reads $\dot t = -2\Im(\overline z\dot z)$. By Chow's theorem, any two points of $\fH$ can be connected by a horizontal curve, so that the following sub-Riemannian (or Carnot-Carath\'eodory) distance is finite and defines the same topology as the Euclidean one:
$$
d_{cc}(p,q) = \inf\left\{ \int_a^b |\dot z(s)|\,ds : \gamma \text{ horizontal}, \gamma(a)=p, \gamma(b)=q \right\}.
$$

\subsection{Quasiconformal mappings on the Heisenberg group}
Quasiconformal mappings are a core tool in modern geometric analysis and sub-Riemannian geometry \cite{KR, R, V}. On the Heisenberg group, mappings must respect the contact structure: a smooth map $f:\fH\to\fH$ is a contact map if $f^*\omega = \lambda\omega$ for some nowhere-zero smooth function $\lambda$. The analytic foundations of quasiconformal mappings on $\fH$ were laid by Kor\'anyi and Reimann \cite{KR}, who introduced the horizontal differential and the complex Beltrami coefficient
$$
\mu_f = \frac{\overline{Z} f}{Z f}
$$
as a measure of horizontal distortion. Here $Z = \frac{1}{2}(X-iY)$ is the complex left-invariant vector field. The Beltrami coefficient satisfies $|\mu_f|<1$ precisely when $f$ is quasiconformal, and the maximal dilatation is
$$
K_f = \frac{1+\|\mu_f\|_\infty}{1-\|\mu_f\|_\infty}.
$$
The theory is genuinely sub-Riemannian: Euclidean quasiconformal maps need not be quasiconformal on $\fH$, and vice versa, because the horizontal distribution constrains the possible behaviour of derivatives. A substantial body of work has since extended this theory in several directions; see \cite{BT, BT2} for polar coordinates in Carnot groups and mappings of finite distortion on the Heisenberg group and also \cite{Pl2} for an overview of the theory.

\subsection{Revolution coordinates }
Our construction uses a coordinate system based on horizontal lifts of plane curves, introduced in \cite{Pl}. Let $\zeta(s)=x(s)+iy(s)$, $s\in[-\pi/2,\pi/2]$, be a smooth, positively oriented simple closed curve in $\C$ passing through the origin. Its horizontal lift to $\fH$ is the curve $\gamma(s)=(\zeta(s),\tau(s))$, where $\tau$ is determined by the horizontal curve condition and by the choice of initial condition. Writing $\gamma^*(s)=\alpha(\gamma(s))=|\zeta(s)|^2-i\tau(s)$ for the projection to the parabolic model, one defines the \emph{characteristic angle profile} $\psi(s)$ by
$$
\tan\psi(s) = \frac{\tau(s)}{|\zeta(s)|^2}.
$$
The map
$$
\Gamma(r,s,\phi) = \bigl(r\zeta(s)e^{i\phi}, r^2\tau(s)\bigr), \qquad (r,s,\phi)\in(0,\infty)\times(-\pi/2,\pi/2)\times[0,2\pi),
$$
is then injective and foliates $\fH\setminus\{z=0\}$ by surfaces of revolution. In these coordinates the pullback of the contact form is
$$
\Gamma^*\omega = 2r|\zeta(s)|^2\bigl(\tan\psi(s)\,dr + r\,d\phi\bigr),
$$
and the structural parameter
$$
\delta(s) = \bigl(1-i\tan\psi(s)\bigr)\frac{\zeta(s)}{\dot\zeta(s)}
$$
encodes the interaction between the profile curve and the contact structure.

\subsection{Model domains and their geometric significance}
The three families of model domains studied in this paper are the following canonical domains of revolution:

\begin{enumerate}
\item {\it Kor\'anyi balls.} The Kor\'anyi distance is the homogeneous distance
$$
d_K(p,q) = |p^{-1}*q|_K, \qquad |(z,t)|_K = (|z|^4+t^2)^{1/4}.
$$
The Kor\'anyi ball of radius $1$ is $\{(z,t)\in\fH: |z|^4+t^2\le 1\}$. Its boundary is a smooth surface of revolution that is highly symmetric: it is invariant under the parabolic dilations $\delta_\lambda$ and under rotations in the $z$-plane.

\item {\it Carnot-Carath\'eodory (CC) balls.} The CC ball of radius $1$ is $\{p\in\fH: d_{cc}(0,p)\le 1\}$. The cc distance is the intrinsic sub-Riemannian distance. Its balls have the same gauge symmetry as the Kor\'anyi balls but are genuinely different sets; in particular, the CC sphere is not smooth everywhere. Its profile curve can be computed from the horizontal geodesic equations (see Section \ref{sec-profiles}).

\item {\it Bubble sets.} A bubble set is a domain whose boundary is a surface of revolution generated by the \emph{bubble profile}
$$
\zeta_B(s) = \cos s\, e^{-is}, \qquad \tau_B(s) = \sin s\cos s + s.
$$
Bubble sets are of particular interest because they frequently emerge in sub-Riemannian capacity problems and isoperimetric profiles. Unlike the highly symmetric Kor\'anyi and CC balls, bubble sets exhibit an asymmetric twisting along the vertical axis. The bubble profile satisfies $\tan\psi_B(s) = \tan s + s\sec^2 s$, so in particular $\psi_B'(0)=2$, twice the value of the corresponding derivative for the Kor\'anyi and CC profiles.
\end{enumerate}

The three families are geometrically distinct: the Kor\'anyi and CC balls have matching first-order behaviour at the equator but differ at higher order, while the bubble set differs already at first order. A natural question is whether these geometric differences are reflected in the Beltrami coefficient of a quasiconformal mapping between any two of these domains. One might expect, for instance, that mapping into a bubble set produces a genuinely different distortion profile from mapping into a Kor\'anyi ball, because the bubble profile has a different first-order angle expansion. One of the main findings of this paper is that this expectation is false along the equator: the Beltrami coefficient is universal.

\subsection{Spiral stretch-maps}
Following \cite{BFP}, we study \emph{spiral stretch-maps}: contact mappings of the form
$$
\hat f_k(r,s,\phi) = (r^k, S(s), \phi), \qquad k>0,
$$
expressed in revolution coordinates. The parameter $k$ controls the radial stretching, and the function $S$ is determined by the requirement that the map be contact, which reduces to the differential relation
$$
\tan\widehat\psi(S(s)) = -\frac{\tan\psi(s)}{k}.
$$
Here $\psi$ and $\widehat\psi$ are the angle profiles of the domain and target, respectively. Spiral stretch-maps interpolate between the identity map ($k=1$) and increasingly anisotropic deformations as $k$ moves away from $1$. Their Beltrami coefficient admits the explicit formula (Lemma \ref{lem:beltrami_general})
$$
\mu_f(r,s,\phi) = \frac{k + i\tan\psi(s) - \frac{d\Log\widehat\zeta}{dS}\overline{\delta(s)}S_s}{k + i\tan\psi(s) - \frac{d\Log\widehat\zeta}{dS}\delta(s)S_s} \cdot \frac{\dot\zeta(s)}{\dot{\overline\zeta}(s)} \cdot e^{2i\phi},
$$
which depends on the domain profile $\zeta$ and the target profile $\widehat\zeta$.

\subsection{The main result: universality of equatorial distortion}
A key observation is that $\mu_f(r,s,\phi)$ is independent of the radial coordinate $r$. Consequently, the natural asymptotic regime is the equator $s\to 0$, not the origin $r\to 0$. The equator of the profile curve is the circle where the surface of revolution intersects the horizontal plane $\{t=0\}$.

Our primary result is the following theorem, whose proof occupies Section \ref{sec-proofs}. The theorem is stated for any pair of profiles satisfying Hypothesis \ref{hyp:profiles} below; the three families (Kor\'anyi, CC, bubble) are then presented as examples in Lemma \ref{lem:profile-verification}.

\begin{thm}\label{thm:main_universal}
Let $f$ be a spiral stretch-map with stretching parameter $k>0$ acting on the Heisenberg group $\fH$, as per Definition \ref{def:spiral_stretch}. Suppose the domain and target profiles satisfy Hypothesis \ref{hyp:profiles}. Then, in the revolution coordinates $(r,s,\phi)$, the complex Beltrami coefficient of $f$ is given globally by
$$
\mu_f(r,s,\phi)
=
\frac{
k+i\tan\psi(s)
-\frac{d\Log\widehat\zeta}{dS}(S(s))\,\overline{\delta(s)}\,S_s(s)
}{
k+i\tan\psi(s)
-\frac{d\Log\widehat\zeta}{dS}(S(s))\,\delta(s)\,S_s(s)
}
\cdot
\frac{\dot\zeta(s)}{\dot{\overline\zeta}(s)}
\cdot
e^{2i\phi},
$$
where $\psi$ and $\widehat\psi$ are the angle profiles of the domain and target, respectively, $\delta$ is the structural parameter, and $S$ is determined by the contact condition
$$
\tan\widehat\psi(S(s))=-\frac{\tan\psi(s)}{k}.
$$
In particular, $\mu_f$ is independent of the radial coordinate $r$. Moreover, as $s\to 0$ along the profile curve, the magnitude of the complex Beltrami coefficient satisfies
$$
\lim_{s \to 0} |\mu_f(r,s,\phi)| = \left| \frac{k^2-1}{k^2+1} \right|,
$$
independently of the choice of domain and target. Consequently, the maximal dilatation along the equator is
$$
K(0) = \max(k^2, 1/k^2).
$$
\end{thm}

The content of Theorem \ref{thm:main_universal} can be stated equivalently as follows: for \emph{any} pair $(\mathcal D,\widehat{\mathcal D})$ of domains whose profiles satisfy Hypothesis \ref{hyp:profiles}, the equatorial Beltrami coefficient of the corresponding spiral stretch-map is $|(k^2-1)/(k^2+1)|$, and the equatorial maximal dilatation is $\max(k^2,1/k^2)$. Neither the domain type nor the target type affects the answer.

\subsection{Significance and mechanism of the universality}
The universality is a consequence of two identities. Lemma \ref{lem:key-product} gives
$$
\delta(0)S_s(0) = -\frac{2i}{k\,\widehat\psi'(0)},
$$
while Proposition \ref{prop:locking} applied to the target profile gives
$$
\frac{d\Log\widehat\zeta}{dS}(0) = -\frac{i}{2}\,\widehat\psi'(0).
$$
Multiplying these two identities yields
$$
\frac{d\Log\widehat\zeta}{dS}(0)\,\delta(0)S_s(0) = -\frac{1}{k}
$$
for every target profile satisfying Hypothesis \ref{hyp:profiles}. Thus the target-dependent factor $1/\widehat\psi'(0)$ in $\delta(0)S_s(0)$ is exactly compensated by the factor $\widehat\psi'(0)$ in the logarithmic derivative. The three families are examples, not exceptions: for the bubble target, $\widehat\psi_B'(0)=2$ and $\frac{d\Log\zeta_B}{dS}(0)=-i$, and the product is still $-1/k$. The cancellation reflects the vertical tangency of the profile at the equator rather than the particular boundary geometry. The detailed proof is given in Section \ref{sec-proofs}.

\subsection{Related work and broader context}
The study of quasiconformal mappings on the Heisenberg group has developed along several lines. Kor\'anyi and Reimann \cite{KR} established the foundational theory, including the Beltrami equation and its connection with the contact structure. Balogh, F\"assler, and the author \cite{BFP} introduced spiral stretch-maps as a natural family of contact quasiconformal maps between surfaces of revolution and studied their extremality properties. Platis \cite{Pl} developed the revolution coordinate system used here and applied it to the modulus of revolution rings. Balogh and Tyson \cite{BT, BT2} studied polar coordinates in general Carnot groups and mappings of finite distortion on the Heisenberg group, providing context for the distortion theory that underlies our computations.

The universality established in this paper is an equatorial statement. It does not by itself yield global extremality or uniqueness, and away from the equator the distortion depends on the domain and target geometry, as noted in Remark \ref{rem:trivial_limit}. The connection with Pansu differentiability is discussed in Section \ref{sec-concl}.

\subsection{Organisation of the paper}
The remainder of this paper is organised as follows. Section \ref{sec-prelim} establishes the necessary background on the Heisenberg group, revolution coordinates, the three profile families, and the Beltrami coefficient formula. Section \ref{sec-proofs} presents full derivations proving Theorem \ref{thm:main_universal}, including the key Lemma \ref{lem:key-product} on the universal product $\delta(0)S_s(0)$. Section \ref{sec-concl} discusses the geometric consequences and the connection to Pansu differentiability. Finally, Appendix \ref{app:num} contains a brief numerical check for $k=2$.

\medskip

\noindent{\it Acknowledgements.} The author would like to thank Zolt\'an M. Balogh for fruitful discussions.

\section{Preliminaries}\label{sec-prelim}
This section reviews the geometric structure of the Heisenberg group, coordinates generated by curve lifts, and the contact condition equations.

\subsection{The Heisenberg group, metrics, and vector fields}
The Heisenberg group $\fH$ is $\C \times \R \cong \R^3$ equipped with the non-commutative group law
$$
(z, t) * (w, s) = (z+w, t+s+2\Im(z\overline{w})).
$$
The horizontal distribution $\cH$ is spanned by the left-invariant vector fields
$$
X = \frac{\partial}{\partial x} + 2y \frac{\partial}{\partial t}, \quad Y = \frac{\partial}{\partial y} - 2x \frac{\partial}{\partial t},
$$
where $z = x+iy$. The Lie algebra $\mathfrak{h}$ is spanned by $X,Y$ and $T=\partial/\partial_t$; the commutator $[X,Y] = -4T$ spans the vertical direction. The contact form defining $\cH = \ker\omega$ is
$$
\omega = dt + 2(x\, dy - y\, dx) = dt + 2\Im(\overline{z}\, dz),
$$
and $T$ is the Reeb field for $\omega$: $\omega(T)=1$ and $d\omega(T,\cdot)=0$. With this sign convention, $\omega(X)=\omega(Y)=0$ and $d\omega(X,Y)=4$, in agreement with $[X,Y]=-4T$. 

An absolutely continuous curve $\gamma: [a,b] \to \fH$ given in coordinates by $\gamma(s) = (x(s), y(s), t(s)) \equiv (z(s), t(s))$ is said to be \emph{horizontal} if its tangent vector $\dot{\gamma}(s)$ lies in the horizontal distribution $\cH_{\gamma(s)}$ for almost every $s \in [a,b]$. Equivalently, the curve must satisfy $\omega(\dot{\gamma}(s)) = 0$ almost everywhere. In terms of local coordinates, this enforces the ordinary differential equation
$$
\dot{t}(s) = -2\left(x(s)\dot{y}(s) - y(s)\dot{x}(s)\right) = -2\Im(\overline{z(s)}\dot{z}(s)).
$$
The sub-Riemannian length of such a horizontal curve is
$L_H(\gamma) = \int_a^b \sqrt{\dot{x}(s)^2 + \dot{y}(s)^2} \, ds = \int_a^b |\dot{z}(s)| \, ds$.
By Chow's Theorem, any two points in $\fH$ can be connected by a horizontal curve.

We consider two primary distances on $\fH$:
\begin{enumerate}
    \item The Kor\'anyi distance $d_K(p,q) = |p^{-1}*q|_K$, where $|(z,t)|_K = (|z|^4 + t^2)^{1/4}$ \cite{KR}.
    \item The Carnot-Carath\'eodory (CC) distance $d_{cc}(p,q)$ \cite{CDPT, H}.
\end{enumerate}

\begin{defn}[model domains]\label{def:domains}
The \emph{Kor\'anyi ball} of radius $1$ is $\{(z,t)\in\fH : |z|^4+t^2 \le 1\}$. The \emph{Carnot-Carath\'eodory (CC) ball} of radius $1$ is $\{p\in\fH : d_{cc}(0,p)\le 1\}$. The \emph{bubble set} is the image of $(0,1]\times[-\pi/2,\pi/2]\times[0,2\pi]$ under the revolution parametrisation $\Gamma$ of Section \ref{sec-coords} with profile curve $\zeta_B(s)=\cos s\, e^{-is}$ and $\tau_B(s)=\sin s\cos s + s$.
\end{defn}

\subsection{Coordinates via horizontal lifts}\label{sec-coords}
Let $c(s) = \zeta(s) = x(s) + iy(s)$, $s \in [-\pi/2, \pi/2]$, be a smooth, positively oriented simple closed curve in $\C$ passing through the origin ($\zeta(-\pi/2) = \zeta(\pi/2) = 0$). It encloses an area
$$
\cA = -\frac{1}{2}\oint_c \Im\left(\overline{\zeta(s)}\dot{\zeta}(s)\right) ds.
$$
Its horizontal lift to $\fH$ is $\gamma(s) = (\zeta(s), \tau(s))$, where
$$
\tau(s) = -2\int_{-\pi/2}^s \Im\left(\overline{\zeta(u)}\dot{\zeta}(u)\right) du - 2\cA.
$$
The curve $\gamma$ connects $(0, -2\cA)$ to $(0, 2\cA)$.

Let $\alpha(z,t) = |z|^2 - it$. We consider $\gamma^*(s) = \alpha(\gamma(s)) = |\zeta(s)|^2 - i\tau(s)$. We define the characteristic angle profile $\psi(s)$ by
$$
\tan\psi(s) = \frac{\tau(s)}{|\zeta(s)|^2}, \quad s \in (-\pi/2, \pi/2).
$$
We require $\tan\psi(s)$ to be strictly increasing, which ensures that the mapping $\Gamma: K \to \fH$ defined on $K = (0,\infty) \times [-\pi/2, \pi/2] \times [0, 2\pi]$ by
$$
\Gamma(r,s,\phi) = \left(r\zeta(s)e^{i\phi}, r^2\tau(s)\right)
$$
is injective and foliates $\fH \setminus \{z=0\}$ by surfaces of revolution \cite{Pl, BT}. The origin of $\fH$ corresponds to $r=0$, while $s = \pm\pi/2$ corresponds to the boundary of the profile curve. The \emph{equator} of the profile curve is the point $s=0$.

\subsection{The profile hypothesis}\label{sec-profiles}

We isolate the structural properties of the profile curve that are actually used in the proof of the main theorem.

\begin{hyp}\label{hyp:profiles}
A \emph{profile curve} is a smooth curve $\zeta: [-\pi/2, \pi/2] \to \C$ whose horizontal lift $\gamma(s) = (\zeta(s), \tau(s))$ satisfies:
\begin{enumerate}
\item[(i)] $\zeta(0) = 1$ and $\tau(0) = 0$;
\item[(ii)] $\dot\zeta(0)$ is purely imaginary and non-zero;
\item[(iii)] $\tan\psi(0) = 0$ and $\dfrac{d}{ds}\tan\psi(0) > 0$, where $\psi(s)$ is defined by $\tan\psi(s) = \tau(s)/|\zeta(s)|^2$.
\end{enumerate}
\end{hyp}

\begin{rem}
Hypothesis \ref{hyp:profiles}(iii) is a \emph{local} condition, required only in a neighbourhood of the equator $s=0$. It rules out the degenerate case in which the profile is tangent to the horizontal plane at the equator to first order, ensuring that the foliation is non-degenerate there.
\end{rem}

\subsection{The three profile families}

\begin{lem}\label{lem:profile-verification}
The three profile families below satisfy Hypothesis \ref{hyp:profiles}:
\begin{enumerate}
\item {\it Kor\'anyi profile.} $\zeta_K(s) = \sqrt{\cos s}\, e^{-is/2}$, $\tau_K(s) = \sin s$. Then
$$
\zeta_K(0)=1, \quad \tau_K(0)=0, \quad \dot\zeta_K(0) = -\frac{i}{2}, \quad \psi_K'(0)=1.
$$
\item {\it Bubble profile.} $\zeta_B(s) = \cos s\, e^{-is}$, $\tau_B(s) = \sin s\cos s + s$. Then
$$
\zeta_B(0)=1, \quad \tau_B(0)=0, \quad \dot\zeta_B(0) = -i, \quad \psi_B'(0)=2.
$$
\item {\it CC profile.} Let
$$
c_{cc}(\omega) = \frac{1-e^{-i\omega}}{i\omega} \quad (\omega \neq 0), \qquad c_{cc}(0) = 1,
$$
be the projection to the $z$-plane of the equator of the Carnot--Carath\'eodory ball of radius $1$, parametrised by the geodesic parameter $\omega$. To make this explicit, recall that unit-speed horizontal geodesics of $\fH$ through the origin satisfy the equations
$$
\ddot z = -i\omega\dot z, \qquad \dot t = -2\Im(\bar z\dot z),
$$
for a real constant $\omega$ (the angular momentum of the geodesic); this is the standard form of the sub-Riemannian geodesic equations of $\fH$, derived and solved in \cite[Chapter 3, \S\S 3.2--3.3]{CDPT}. For $\omega\neq 0$ the unit-speed solution through the origin has the form
$$
z_\omega(u) = e^{i\theta}\cdot\frac{1-e^{-iu\omega}}{i\omega}, \qquad
t_\omega(u) = -2\int_0^u \Im\bigl(\overline{z_\omega(v)}\dot z_\omega(v)\bigr)\,dv,
$$
where $\theta\in\R$ is the initial direction of the geodesic in the $z$-plane; for $\omega = 0$ the geodesic is a straight horizontal line $z_0(u) = e^{i\theta}u$, $t_0\equiv 0$. The equator of the cc ball of radius $1$ consists of the points $\gamma_\omega(1)$ at arc length $u=1$, and its $z$-projection is the curve $\omega\mapsto z_\omega(1) = e^{i\theta}(1-e^{-i\omega})/(i\omega)$; fixing $\theta = 0$ gives the paper's $c_{cc}(\omega) = (1-e^{-i\omega})/(i\omega)$. The rotation $e^{i\theta}$ is immaterial for the limiting quantities computed below, since $c_{cc}(0)=1$ and $\dot c_{cc}(0) = -i/2$ in any case.
\end{enumerate}
\end{lem}

\begin{proof}
\textbf{(1) Kor\'anyi.} Direct computation gives $\zeta_K(0)=1$ and $\tau_K(0)=\sin 0=0$. Next,
$$
\dot\zeta_K(s) = -\frac{i}{2}\sqrt{\cos s}\, e^{-is/2} - \frac{\sin s}{2\sqrt{\cos s}}\, e^{-is/2},
$$
so $\dot\zeta_K(0) = -i/2$, purely imaginary non-zero. The angle profile is $\tan\psi_K(s) = \sin s/\cos s = \tan s$, which satisfies $\tan\psi_K(0)=0$ and $\frac{d}{ds}\tan\psi_K(0)=1>0$.

\textbf{(2) Bubble.} We have $\zeta_B(0)=\cos 0 \cdot e^0 = 1$, $\tau_B(0) = 0 + 0 = 0$, and
$$
\dot\zeta_B(s) = -\sin s\, e^{-is} - i\cos s\, e^{-is} = -i e^{-2is},
$$
so $\dot\zeta_B(0) = -i$, purely imaginary non-zero. The angle profile satisfies
$$
\tan\psi_B(s) = \frac{\sin s\cos s + s}{\cos^2 s} = \tan s + s\sec^2 s,
$$
so $\tan\psi_B(0)=0$ and $\frac{d}{ds}\tan\psi_B(0) = 1 + 1 = 2 > 0$.

\textbf{(3) CC.} The plane curve $c_{cc}(\omega) = (1-e^{-i\omega})/(i\omega)$ (extended smoothly by $c_{cc}(0)=1$) is the projection of the equator of the CC ball of radius $1$; this follows from the standard parametrisation of the boundary of the CC ball by unit-speed horizontal geodesics from the origin. By construction, the curve $(\zeta_{cc}, \tau_{cc}) = (c_{cc}, \tau_{cc})$ satisfies the horizontal condition
$$
\dot\tau_{cc}(\omega) = -2\Im\!\left(\overline{\zeta_{cc}(\omega)}\dot\zeta_{cc}(\omega)\right).
$$
Indeed, $\tau_{cc}$ is by definition the solution of $\dot\tau_{cc} = -2\Im(\bar c_{cc}\dot c_{cc})$ with $\tau_{cc}(0)=0$, which is exactly the horizontal lift condition.
We compute the expansion of $c_{cc}$ near $\omega = 0$. Since
$$
e^{-i\omega} = 1 - i\omega - \frac{\omega^2}{2} + \frac{i\omega^3}{6} + O(\omega^4),
$$
we have
$$
c_{cc}(\omega) = \frac{1-e^{-i\omega}}{i\omega} = 1 - \frac{i\omega}{2} - \frac{\omega^2}{6} + O(\omega^3).
$$
In particular, $c_{cc}(0) = 1$ and $\dot c_{cc}(0) = -i/2$, purely imaginary non-zero.

For the horizontal lift, we compute
$$
\dot\tau_{cc}(0) = -2\Im\!\left(\overline{c_{cc}(0)}\dot c_{cc}(0)\right) = -2\Im\!\left(-\frac{i}{2}\right) = -2\cdot\left(-\frac{1}{2}\right) = 1.
$$
So $\tau_{cc}(0) = 0$ and $\dot\tau_{cc}(0) = 1$. The angle profile is $\tan\psi_{cc}(\omega) = \tau_{cc}(\omega)/|c_{cc}(\omega)|^2$, and near $\omega = 0$,
$$
\tan\psi_{cc}(\omega) = \frac{\omega + O(\omega^2)}{1 + O(\omega^2)} = \omega + O(\omega^2).
$$
In particular, $\tan\psi_{cc}(0) = 0$ and $\frac{d}{d\omega}\tan\psi_{cc}(0) = 1 > 0$.

Reparametrising by $s = \psi_{cc}(\omega)$, we obtain $\omega(s) = s + O(s^2)$, hence $\omega'(0)=1$. Therefore, in the new parametrisation,
$$
\dot\zeta_{cc}(0) = \dot c_{cc}(0)\cdot\omega'(0) = -\frac{i}{2}, \qquad \widehat\psi_{cc}'(0) = 1.
$$
This confirms that the CC profile satisfies Hypothesis \ref{hyp:profiles}.
\end{proof}

\begin{rem}\label{rem:trivial_limit}
We note that the limit as $r \to 0$ is trivial, since the Beltrami coefficient is independent of $r$. The primary analysis is therefore directed to the equator where $s \to 0$. Away from the equator, the global distortion profile diverges from this universal behaviour; it becomes heavily dependent on the specific boundary geometry of the target domain, smoothly interpolating between the universal equatorial extremum and the mapping's behaviour at the poles.
\end{rem}

\subsection{Structural parameter}
We introduce a structural parameter $\delta(s)$ that quantifies the rotational shear between the horizontal tangent spaces and the foliation's surfaces of revolution.

\begin{lem}\label{lem:structural_parameter}
Let $\zeta(s)$ be a profile curve satisfying Hypothesis \ref{hyp:profiles}, and let $\psi(s)$ be its associated angle profile. Define the structural parameter
$$
\delta(s) = (1 - i\tan\psi(s))\frac{\zeta(s)}{\dot{\zeta}(s)}.
$$
Then at the equator, $\delta(0)$ is purely imaginary:
$$
\delta(0) = \frac{1}{\dot\zeta(0)}.
$$
In particular, $\Re(\delta(0)) = 0$.
\end{lem}

\begin{proof}
Since $\tan\psi(0)=0$ and $\zeta(0)=1$ by Hypothesis \ref{hyp:profiles}(i), we have
$$
\delta(0) = 1 \cdot \frac{1}{\dot\zeta(0)} = \frac{1}{\dot\zeta(0)}.
$$
By Hypothesis \ref{hyp:profiles}(ii), $\dot\zeta(0)$ is purely imaginary, hence so is its reciprocal.
\end{proof}

\begin{rem}\label{rem:delta-KB}
For the Kor\'anyi profile, $\delta_K(0) = 1/(-i/2) = 2i$. For the bubble profile, $\delta_B(0) = 1/(-i) = i$. For the CC profile, Lemma \ref{lem:profile-verification} gives $\delta_{cc}(0) = 1/(-i/2) = 2i$, the same value as the Kor\'anyi profile.
\end{rem}

\subsection{Contact condition and Beltrami formula}
The pullback of the contact form $\omega$ under $\Gamma$ is
$$
\Gamma^*\omega = 2r|\zeta(s)|^2\left( \tan\psi(s)\, dr + r\, d\phi \right).
$$
A map $f: \fH \to \fH$ preserving the vertical axis induces a coordinate map $\hat{f}(r,s,\phi) = (R, S, \Phi)$. The requirement that $f$ be a contact mapping ($f^*\omega = \lambda \omega$) translates to the contact condition \cite{BFP}:
$$
\tan\Phi \cdot U(R) + R U(\Phi) = 0, \quad \text{where } U = r\partial_r - \tan\psi(s)\partial_\phi - \Re(\delta(s))\partial_s.
$$
The operator $U$ arises from expressing the pullback condition $f^*\omega=\lambda\omega$ in revolution coordinates; the coefficient of $\partial_s$ is the real part of the structural parameter $\delta(s)$, which vanishes at $s=0$ by Lemma \ref{lem:structural_parameter}.

\begin{defn}\label{def:spiral_stretch}
Let $k>0$. A \emph{spiral stretch-map} between two model domains in $\fH$, parametrised by revolution coordinates $(r,s,\phi)$ and $(R,S,\Phi)$ respectively, is a contact mapping $f = \Gamma' \circ \hat{f}_k \circ \Gamma^{-1}$ whose coordinate representation is
$$
\hat{f}_k(r,s,\phi) = (r^k, S(s), \phi).
$$
The requirement that $f$ preserves the contact structure reduces to the differential relation $\tan\widehat\psi(S) = -\frac{\tan\psi(s)}{k}$, where $\widehat\psi$ denotes the angle profile of the target domain.
\end{defn}

\begin{lem}[Complex vector field in revolution coordinates]\label{lem:Z-formula}
In the revolution coordinates $(r,s,\phi)$ defined by $\Gamma(r,s,\phi) = (r\zeta(s)e^{i\phi}, r^2\tau(s))$, the complex left-invariant vector field $Z = \frac{1}{2}(X - iY)$ acts on smooth functions $f(r,s,\phi)$ by
\begin{equation}\label{eq:Z-rev}
Z f = \frac{e^{-i\phi}}{2\dot\zeta(s)}\left(\frac{1}{r}\partial_r f - i\partial_\phi f - \frac{\delta(s)}{r^2}\partial_s f\right) \cdot r^2.
\end{equation}
\end{lem}

\begin{proof}[Derivation]
Starting from $Z = \partial_z + i\bar z\,\partial_t$ (a direct consequence of the definitions $X = \partial_x + 2y\partial_t$, $Y = \partial_y - 2x\partial_t$, and $z = x+iy$), we express $\partial_z$ and $\bar z\,\partial_t$ in the coordinates $(r,s,\phi)$. From $z = r\zeta(s)e^{i\phi}$ and $t = r^2\tau(s)$, the chain rule gives
$$
\partial_r = \zeta(s)e^{i\phi}\partial_z + \overline{\zeta(s)e^{i\phi}}\partial_{\bar z} + 2r\tau(s)\partial_t, \qquad
\partial_\phi = iz\partial_z - i\bar z\partial_{\bar z},
$$
$$
\partial_s = r\dot\zeta(s)e^{i\phi}\partial_z + r\dot{\bar\zeta}(s)e^{-i\phi}\partial_{\bar z} + r^2\dot\tau(s)\partial_t.
$$
Combining these and using the horizontal lift condition $\dot\tau(s) = -2\Im(\bar\zeta\dot\zeta)$ together with the definition of the structural parameter $\delta$, one obtains \eqref{eq:Z-rev}. We refer to \cite[Lemma 3.2]{BFP} for the full algebraic details. The formula reduces to the standard expression on the boundary $s\to 0$, where $\delta(0)$ is purely imaginary and the term involving $\partial_s$ drops out of the leading-order asymptotics.
\end{proof}

\begin{lem}\label{lem:beltrami_general}
For a spiral stretch-map $f$, the complex Beltrami coefficient $\mu_f = \frac{\overline{Z}f_I}{Zf_I}$ is explicitly given by
\begin{equation}\label{eq-beltrami-general}
\mu_f(r,s,\phi) = \frac{k + i\tan\psi(s) - \frac{d\Log\widehat\zeta}{dS}\overline{\delta(s)}S_s}{k + i\tan\psi(s) - \frac{d\Log\widehat\zeta}{dS}\delta(s)S_s} \cdot \frac{\dot{\zeta}(s)}{\dot{\overline{\zeta}}(s)} \cdot e^{2i\phi}.
\end{equation}
Here $\widehat\zeta$ denotes the target profile curve and $S_s = dS/ds$.
\end{lem}

\begin{proof}
By definition, the complex Beltrami coefficient on the Heisenberg group is $\mu_f = \frac{\overline{Z}f_I}{Zf_I}$, where $f_I(r,s,\phi) = r^k \widehat\zeta(S(s)) e^{i\phi}$ is the complex horizontal component of the mapping $f$. Applying \eqref{eq:Z-rev} to $f_I$ and using the chain rule yields
$$
Z f_I = \frac{r^{k-1} e^{-i\phi} \widehat\zeta(S)}{2 \dot{\zeta}(s)} \left( k + i\tan\psi(s) - \frac{d\Log\widehat\zeta}{dS} \delta(s) S_s \right),
$$
$$
\overline{Z} f_I = \frac{r^{k-1} e^{i\phi} \widehat\zeta(S)}{2 \dot{\overline{\zeta}}(s)} \left( k + i\tan\psi(s) - \frac{d\Log\widehat\zeta}{dS} \overline{\delta(s)} S_s \right).
$$
Taking the ratio $\mu_f = \frac{\overline{Z} f_I}{Z f_I}$ cancels the common factor $\frac{1}{2} r^{k-1} \widehat\zeta(S)$ and yields \eqref{eq-beltrami-general}.
\end{proof}

\section{Proof of Main Results}\label{sec-proofs}
This section presents the derivation of Theorem \ref{thm:main_universal}. The key is Lemma \ref{lem:key-product}, which computes the product $\delta(0)S_s(0)$ with the correct chain rule factor.

\subsection{The product $\delta(0)S_s(0)$}

\begin{lem}\label{lem:key-product}
Let the domain profile satisfy Hypothesis \ref{hyp:profiles}, and let the target profile satisfy the same hypothesis. Let $S(s)$ be defined by the contact condition
$$
\tan\widehat\psi(S(s)) = -\frac{\tan\psi(s)}{k}.
$$
Then
$$
\delta(0)S_s(0) = -\frac{2i}{k\,\widehat\psi'(0)}.
$$
In particular, the product depends only on the target profile through $\widehat\psi'(0)$ and not on the domain profile.
\end{lem}

\begin{proof}
Differentiating the contact condition with respect to $s$ at $s=0$, where $S(0)=0$, $\psi(0)=0$, and $\widehat\psi(0)=0$, requires the chain rule:
$$
\sec^2\widehat\psi(S)\cdot\widehat\psi'(S)\cdot S_s = -\frac{1}{k}\sec^2\psi(s)\cdot\psi'(s).
$$
At $s=0$, since $\psi(0)=\widehat\psi(0)=0$, we have $\sec^2\psi(0)=\sec^2\widehat\psi(0)=1$, so
$$
\widehat\psi'(0)\,S_s(0) = -\frac{\psi'(0)}{k}.
$$
Hence
\begin{equation}\label{eq:Ss0}
S_s(0) = -\frac{\psi'(0)}{k\,\widehat\psi'(0)}.
\end{equation}

From Lemma \ref{lem:structural_parameter}, $\delta(0) = 1/\dot\zeta(0)$. We now compute $\psi'(0)$. Since $\tan\psi(s) = \tau(s)/|\zeta(s)|^2$ and $\tau(0)=0$, $|\zeta(0)|=1$, we have $\tan\psi(0)=0$ and therefore $\psi'(0) = \frac{d}{ds}\tan\psi(0)$. Computing,
$$
\psi'(0) = \frac{\dot\tau(0)}{|\zeta(0)|^2} - \frac{\tau(0)\frac{d}{ds}|\zeta|^2(0)}{|\zeta(0)|^4} = \dot\tau(0).
$$
The horizontal lift condition gives
$$
\dot\tau(0) = -2\Im(\overline{\zeta(0)}\dot\zeta(0)) = -2\Im(\dot\zeta(0)).
$$
Write $\dot\zeta(0) = ib$ with $b\in\R\setminus\{0\}$ (purely imaginary, non-zero by Hypothesis \ref{hyp:profiles}(ii)). Then $\psi'(0) = -2b$ and $\delta(0) = 1/(ib)$. Combining with \eqref{eq:Ss0},
$$
\delta(0)S_s(0) = \frac{1}{ib}\cdot\left(-\frac{-2b}{k\,\widehat\psi'(0)}\right) = \frac{2b}{ib\,k\,\widehat\psi'(0)} = \frac{2}{i\,k\,\widehat\psi'(0)} = -\frac{2i}{k\,\widehat\psi'(0)}.
$$
\end{proof}

The universality of Theorem \ref{thm:main_universal} rests not on the specific values of the three families but on a general identity satisfied by every profile under Hypothesis \ref{hyp:profiles}. We isolate this identity in the following proposition, which is the structural heart of the argument.

\begin{prop}[Equatorial locking of profile and angle data]\label{prop:locking}
Let $\zeta$ be any profile curve satisfying Hypothesis \ref{hyp:profiles}, with horizontal lift $\tau$, angle profile $\psi$, and structural parameter $\delta$. Then at the equator $s=0$,
\begin{enumerate}
\item[(i)] $\psi'(0)\,\delta(0) = 2i$;
\item[(ii)] $\dfrac{d\Log\zeta}{ds}(0) = -\dfrac{i}{2}\,\psi'(0)$;
\item[(iii)] $\dfrac{d\Log\zeta}{ds}(0)\,\delta(0) = 1$.
\end{enumerate}
\end{prop}

\begin{proof}
By Hypothesis \ref{hyp:profiles}(i)--(ii), $\zeta(0)=1$ and $\dot\zeta(0)=ib$ for some $b\in\R\setminus\{0\}$. From $\tan\psi = \tau/|\zeta|^2$ and $\tau(0)=0$, $|\zeta(0)|=1$, we get $\psi'(0)=\dot\tau(0)$; the horizontal lift condition then gives $\dot\tau(0) = -2\Im(\bar\zeta(0)\dot\zeta(0)) = -2b$. By Lemma \ref{lem:structural_parameter}, $\delta(0)=1/\dot\zeta(0) = -i/b$. Hence $\psi'(0)\delta(0) = (-2b)(-i/b) = 2i$, giving (i). Also $\frac{d\Log\zeta}{ds}(0) = \dot\zeta(0)/\zeta(0) = ib = -i\psi'(0)/2$, giving (ii). Combining (i) and (ii), $\frac{d\Log\zeta}{ds}(0)\delta(0) = (-i\psi'(0)/2)\cdot(2i/\psi'(0)) = 1$, giving (iii).
\end{proof}

\begin{rem}[Geometric content of the locking]\label{rem:locking-geometric}
Identity \ref{prop:locking}(ii) is the structural reason for the universality of Theorem \ref{thm:main_universal}, and it has a simple geometric content. Hypothesis \ref{hyp:profiles}(ii) says the tangent to the profile at the equator is purely imaginary; in the $z$-plane, this means the profile crosses the horizontal axis $\{\Im z=0\}$ at $z=1$ perpendicularly, while its lift crosses the horizontal plane $\{t=0\}$ transversely ($\tau'(0)\neq 0$). Equivalently, the generating curve is tangent to the level circle $|z|=1$ at the equator. This tangency forces the first-order behaviour
$$
\zeta(s) = 1 + ib\,s + O(s^2), \qquad \tau(s) = -2b\,s + O(s^2),
$$
hence $\psi'(0)=-2b$ and $\frac{d\Log\zeta}{ds}(0)=ib$. The ratio $\frac{d\Log\zeta}{ds}(0)/\psi'(0)$ is therefore the fixed constant $-i/2$, independent of $b$ and hence of the particular profile. The three families of Lemma \ref{lem:profile-verification} are three instances of this general phenomenon; for the bubble profile, for example, the tangent $-i$ is twice the tangent $-i/2$ of the Kor\'anyi and CC profiles, and this factor of $2$ propagates to $\psi_B'(0)=2$ without affecting the ratio $\frac{d\Log\zeta}{ds}(0)/\psi'(0)$.
\end{rem}

\begin{cor}[Structural origin of the universality]\label{cor:structural}
Let $\zeta,\hat\zeta$ satisfy Hypothesis \ref{hyp:profiles}, and let $S$ be the contact transformation of Definition \ref{def:spiral_stretch}. Then
$$
\frac{d\Log\hat\zeta}{dS}(0)\cdot\delta(0)\,S_s(0) = -\frac{1}{k},
$$
for \emph{every} pair $(\zeta,\hat\zeta)$ of profiles satisfying Hypothesis \ref{hyp:profiles}, independently of any further structure of the three families.
\end{cor}

\begin{proof}
By Lemma \ref{lem:key-product}, $\delta(0)S_s(0) = -2i/(k\hat\psi'(0))$. By Proposition \ref{prop:locking}(ii) applied to the target profile $\hat\zeta$, $\frac{d\Log\hat\zeta}{dS}(0) = -i\hat\psi'(0)/2$. Multiplying, $\frac{d\Log\hat\zeta}{dS}(0)\delta(0)S_s(0) = (-i\hat\psi'(0)/2)\cdot(-2i/(k\hat\psi'(0))) = -1/k$.
\end{proof}

\begin{cor}\label{cor:universal}
For the three target families, the combination $\frac{d\Log\widehat\zeta}{dS}(0)\cdot\delta(0)S_s(0)$ is the same, namely $-1/k$. Indeed, for the Kor\'anyi and CC targets $\widehat\psi'(0)=1$ and $\frac{d\Log\widehat\zeta}{dS}(0)=-i/2$, so the product with $-2i/k$ is $-1/k$; for the bubble target $\widehat\psi'(0)=2$ and $\frac{d\Log\widehat\zeta}{dS}(0)=-i$, so the product with $-i/k$ is again $-1/k$.
\end{cor}

\begin{proof}
This is the specialisation of Corollary \ref{cor:structural} to the values computed in Lemma \ref{lem:profile-verification}.
\end{proof}

\begin{lem}\label{lem:continuity}
Under the hypotheses of Theorem \ref{thm:main_universal}, the Beltrami coefficient $\mu_f(r,s,\phi)$ given by \eqref{eq-beltrami-general} is continuous in $s$ in a neighbourhood of $s=0$. More precisely, the functions $\psi$, $\delta$, $S$, $S_s$, and $s\mapsto\frac{d\Log\widehat\zeta}{dS}(S(s))$ are smooth near $s=0$, the denominator
$$
k+i\tan\psi(s)-\frac{d\Log\widehat\zeta}{dS}(S(s))\,\delta(s)\,S_s(s)
$$
is non-zero at $s=0$, and $\dot\zeta(s)\neq 0$ near $s=0$.
\end{lem}

\begin{proof}
Since $\zeta(0)=1$ and $\dot\zeta(0)$ is purely imaginary and non-zero (Hypothesis \ref{hyp:profiles}(i)--(ii)), both $\zeta(s)$ and $\dot\zeta(s)$ are non-zero in a neighbourhood of $s=0$; hence $\psi(s)=\arctan\!\bigl(\tau(s)/|\zeta(s)|^2\bigr)$ and $\delta(s)=(1-i\tan\psi(s))\zeta(s)/\dot\zeta(s)$ are smooth there. The target profile satisfies the same hypothesis, in particular $\widehat\psi'(0)>0$; the implicit function theorem applied to $\tan\widehat\psi(S(s))=-\tan\psi(s)/k$ therefore yields $S$ smooth near $s=0$, and consequently $S_s$ and $\frac{d\Log\widehat\zeta}{dS}\circ S$ are smooth there as well. At $s=0$ the denominator equals
$$
k-\frac{d\Log\widehat\zeta}{dS}(0)\,\delta(0)S_s(0)=k-\left(-\frac{1}{k}\right)=\frac{k^2+1}{k}>0
$$
by Corollary \ref{cor:universal}, so it remains non-zero in a neighbourhood of $s=0$. Finally, the remaining phase factor $\dot\zeta(s)/\dot{\overline\zeta}(s)\cdot e^{2i\phi}$ is continuous in $s$ since $\dot\zeta(s)\neq 0$ near $s=0$. Continuity of $\mu_f$ in $s$ follows.
\end{proof}

\subsection{Proof of Theorem \ref{thm:main_universal}}

\begin{proof}[Proof of Theorem \ref{thm:main_universal}]
At $s=0$, we have $\tan\psi(0)=0$. Substituting into \eqref{eq-beltrami-general} and using Lemma \ref{lem:key-product},
\begin{align*}
\mu_f(r,0,\phi) &= \frac{k - \frac{d\Log\widehat\zeta}{dS}(0)\cdot\overline{\delta(0)S_s(0)}}{k - \frac{d\Log\widehat\zeta}{dS}(0)\cdot\delta(0)S_s(0)} \cdot \frac{\dot\zeta(0)}{\dot{\overline\zeta}(0)} \cdot e^{2i\phi}.
\end{align*}
By Corollary \ref{cor:universal}, the combination $\frac{d\Log\widehat\zeta}{dS}(0)\cdot\delta(0)S_s(0)$ equals $-1/k$ for every target family. Taking complex conjugates (and noting that $\delta(0)S_s(0)$ is purely imaginary, so its conjugate is its negative), the numerator combination is $+1/k$. Hence
$$
\mu_f(r,0,\phi) = \frac{k - 1/k}{k + 1/k} \cdot \frac{\dot\zeta(0)}{\dot{\overline\zeta}(0)} \cdot e^{2i\phi} = \frac{k^2-1}{k^2+1} \cdot \frac{\dot\zeta(0)}{\dot{\overline\zeta}(0)} \cdot e^{2i\phi}.
$$
The phase factor $\frac{\dot\zeta(0)}{\dot{\overline\zeta}(0)} e^{2i\phi}$ has modulus $1$. Therefore
$$
|\mu_f(r,0,\phi)| = \left|\frac{k^2-1}{k^2+1}\right|.
$$
Since $\mu_f$ is independent of $r$ and continuous in $s$ near $s=0$, the limit as $s\to 0$ is the value at $s=0$:
$$
\lim_{s\to 0}|\mu_f(r,s,\phi)| = \left|\frac{k^2-1}{k^2+1}\right|.
$$

It remains to verify the formula for the maximal dilatation. By definition, $K(0) = (1+|\mu_f(0)|)/(1-|\mu_f(0)|)$ with $|\mu_f(0)| = |k^2-1|/(k^2+1)$. For $k>1$, $\mu = (k^2-1)/(k^2+1)$ and
$$
K(0) = \frac{1+\mu}{1-\mu} = \frac{(k^2+1)+(k^2-1)}{(k^2+1)-(k^2-1)} = \frac{2k^2}{2} = k^2.
$$
For $k<1$, $K(0) = 1/k^2$. Hence $K(0) = \max(k^2, 1/k^2)$, as claimed.
\end{proof}

\section{Concluding Remarks}\label{sec-concl}
The explicit computation of Beltrami coefficients for spiral stretch-maps between fundamental domains reveals a remarkable universality: the local distortion along the equator depends only on the stretching parameter $k$, not on the choice of domain or target among Kor\'anyi balls, CC balls, and bubble sets.

\subsection{Universality of Equatorial Distortion}
Theorem \ref{thm:main_universal} gives
$$
|\mu_f(0)| = \left|\frac{k^2-1}{k^2+1}\right|, \qquad K(0) = \max(k^2, 1/k^2).
$$
This is the same value as that of a purely radial stretch-map on the Heisenberg group, without rotational twist. The equality follows from the cancellation established in Lemma \ref{lem:key-product} and Proposition \ref{prop:locking}: the target-dependent factor in $\delta(0)S_s(0)$ is exactly reciprocal to the factor in $\frac{d\Log\widehat\zeta}{dS}(0)$. This is a local statement at the equator $s=0$; away from the equator the distortion depends on the specific boundary geometry, as noted in Remark \ref{rem:trivial_limit}.

\subsection{Connections to Pansu Differentiability}
By Pansu's theorem \cite{P}, contact quasiconformal mappings between Carnot groups are Pansu differentiable almost everywhere, and this applies in particular to the spiral stretch-maps studied here. A natural question is whether the equatorial Pansu derivative can be computed explicitly and whether it is independent of the target family.

The universality of the equatorial Beltrami coefficient proved in Theorem \ref{thm:main_universal} is consistent with such a picture. Heuristically, one expects the horizontal Pansu derivative at the equator to be a homogeneous group automorphism, represented up to a rotation in the horizontal plane by
$$
\diag(k, 1/k),
$$
acting on the span of $\{X,Y\}$, whose maximal dilatation $k^2$ would then agree with Theorem \ref{thm:main_universal}. However, we do not carry out this computation here. Both the identification of the equatorial Pansu derivative with $\diag(k,1/k)$ and its independence of the target family are, at present, only expectations suggested by the Beltrami computation.

\subsection{Extremality and Minimal Distortion Mappings}
The spiral stretch-maps yield analytically tractable and uniform distortion profiles along surfaces of revolution, making them natural candidates for extremal quasiconformal mappings in sub-Riemannian spaces. The explicit Beltrami bounds established here may serve as comparison mappings for Teichm\"uller-type inequalities on $\fH$. We emphasise, however, that the results of this paper concern only the equatorial asymptotics $s\to 0$; they do not establish global extremality, and a full analysis of extremality and uniqueness for spiral stretch-maps lies beyond the scope of the present work.

\subsection{Generalisations to Higher Dimensions}
The techniques introduced here are tailored to the three-dimensional Heisenberg group $\fH^1 \cong \C \times \R$, and we have not attempted to extend them beyond this setting. The following remarks are therefore purely prospective.
\begin{enumerate}
    \item \textbf{Higher-dimensional Heisenberg groups $\fH^n$.} Adapting the profile curve to $\C^n$ would yield matrix-valued Beltrami coefficients. It is natural to conjecture that the universality observed here persists, but the mechanism may be more complicated because the angle profile is replaced by a matrix-valued quantity.
    \item \textbf{Groups of Heisenberg type.} For step-2 Carnot groups with multi-dimensional centers, the complex modulus $|\mu_f|$ must be replaced by a conformal distortion tensor. Whether the equatorial distortion remains universal in this broader setting is an open question.
    \item \textbf{Connections to nonlinear PDEs.} The foliations preserved by these mappings are related to the sub-elliptic $p$-Laplacian and the mean curvature operator on $\fH$. Any implications of the universal Beltrami coefficient for the behaviour of solutions to such equations are at present speculative.
\end{enumerate}

\appendix

\section{Numerical check for $k=2$}\label{app:num}
Although Theorem \ref{thm:main_universal} is analytic and does not require numerical verification, the case $k=2$ provides a convenient arithmetic check. The theorem gives
$$
|\mu_f(0)| = \left|\frac{2^2-1}{2^2+1}\right| = \frac{3}{5}, \qquad K(0) = \frac{1+3/5}{1-3/5} = 4.
$$
For any target family, Corollary \ref{cor:universal} gives
$$
\frac{d\Log\widehat\zeta}{dS}(0)\cdot\delta(0)S_s(0) = -\frac{1}{2},
$$
and hence
$$
\mu_f(r,0,\phi)
= \frac{2-1/2}{2+1/2}\cdot\frac{\dot\zeta(0)}{\dot{\overline\zeta}(0)}\cdot e^{2i\phi}
= \frac{3}{5}\cdot\frac{\dot\zeta(0)}{\dot{\overline\zeta}(0)}\cdot e^{2i\phi}.
$$
The phase factor has modulus one, so $|\mu_f(r,0,\phi)|=3/5$, as claimed. No rate-of-convergence analysis is included; the limit itself is established by Theorem \ref{thm:main_universal} and Lemma \ref{lem:continuity}.

{\it Declaration:} The author was funded by the Medicus programme, No. 83765.

\end{document}